 \documentclass{amsart}

\usepackage{upref, amssymb, hyperref}

\newcommand{\Siegel}{\mathcal{S}}

\DeclareMathOperator{\GL}{GL}
\newcommand{\opencover}{\mathcal{V}}

\newcommand{\A}{\mathbf{A}}
\newcommand{\G}{\mathbf{G}}
\newcommand{\T}{\mathbf{T}}
\newcommand{\Pg}{\mathbf{P}}

\newcommand{\Lie}[1]{\mathfrak{\lowercase{#1}}}

\newcommand{\closure}[1]{\overline{#1}}
\newcommand{\cover}[1]{\widetilde{#1}}

\newcommand{\integer}{{\mathord{\mathbb{Z}}}}
\newcommand{\rational}{{\mathord{\mathbb{Q}}}}
\newcommand{\real}{{\mathord{\mathbb{R}}}}

\newcommand{\Qrank}{\mathop{\mbox{\upshape$\rational$-rank}}}
\newcommand{\Rrank}{\mathop{\mbox{\upshape$\real$-rank}}}

\newcommand{\bigset}[2]{\left\{\, #1 
 \mathrel{\left| \vphantom {\left\{ #1 \mid #2 \right\} }
 \right.} #2 \,\right\} }

\renewcommand{\see}[1]{{\upshape(}see~\ref{#1}{\upshape)}}

 \newcommand{\seeSect}[1]{{\upshape(}see~\S\ref{#1}{\upshape)}}

\newcommand{\pref}[1]{{\upshape(}\ref{#1}{\upshape)}}
\newcommand{\fullref}[2]{{\ref{#1}\pref{#1-#2}}}

\swapnumbers

\newtheorem{thm}[equation]{Theorem}
\newtheorem{prop}[equation]{Proposition}
\newtheorem{cor}[equation]{Corollary}
\newtheorem{lem}[equation]{Lemma}
\newtheorem{conj}[equation]{Conjecture}
\newtheorem{fact}[equation]{Fact}

\theoremstyle{definition}

\newtheorem{notation}[equation]{Notation}
\newtheorem{rem}[equation]{Remark}

\newtheorem{ack}[equation]{Acknowledgments}

\numberwithin{equation}{section}

\newcommand{\thmrefer}[1]{\renewcommand\theequation
  {\protect\ref{#1}$'$}\addtocounter{equation}{-1}}

 \newcounter{case}

\newenvironment{case}[1][\unskip]{\refstepcounter{case}
 \setcounter{subcase}{0}\em
 \medskip \noindent Case \thecase\ #1.\ }{\unskip\upshape}
 \renewcommand{\thecase}{\arabic{case}}

 \newcounter{subcase}
 
 \renewcommand{\thesubcase}{\arabic{case}.\arabic{subcase}}

\begin{document}

\title[Divergent torus orbits in homogeneous spaces]
 {Divergent torus orbits in \\ homogeneous spaces of $\mathbb{Q}$-rank
two}

\author{Pralay Chatterjee}

\address{Department of Mathematics, Oklahoma State University,
Stillwater, OK 74078}

\address{Mathematics Department, Rice University,
Houston, TX 77005}

\email{pralay@math.rice.edu}

\author{Dave Witte Morris}

\address{Department of Mathematics, Oklahoma State University,
Stillwater, OK 74078}

\address{Department of Mathematics and Computer Science, University of
Lethbridge, Lethbridge, AB T1K~3M4, Canada}

\email{Dave.Morris@uleth.ca}

\urladdr{http://people.uleth.ca/$\sim$dave.morris/}

\date{September 7, 2004. 
 Submitted to \emph{Israel Journal of Mathematics} 
 in September 2004.} 

\begin{abstract}
 Let ${\bf G}$ be a semisimple algebraic $\mathbb{Q}$-group, let $\Gamma$
be an arithmetic subgroup of~${\bf G}$, and let ${\bf T}$ be an
$\mathbb{R}$-split torus in~${\bf G}$. We prove that if there is a
divergent ${\bf T}_{\mathbb{R}}$-orbit in $\Gamma \backslash
{\bf G}_{\mathbb{R}}$, and $\mathop{\text{\rm$\mathbb{Q}$-rank}} {\bf G}
\le 2$, then $\mathop{\rm dim}{\bf T} \le
\mathop{\text{\rm$\mathbb{Q}$-rank}} {\bf G}$. This provides a partial
answer to a question of G.~Tomanov and B.~Weiss.
 \end{abstract}

\maketitle

\section{Introduction}

 Let ${\bf G}$ be a semisimple algebraic $\rational$-group, let $\Gamma$
be an arithmetic subgroup of~$\G$, and let $\T$ be an
$\real$-split torus in~$\G$.  The $\T_{\real}$-orbit of a point
$\Gamma x_0$ in $X = \Gamma \backslash \G_{\real}$ is \emph{divergent} if
the natural orbit map $\T_{\real} \to X \colon t \mapsto \Gamma x_0 t$ is
proper. G.~Tomanov and B.~Weiss \cite[p.~389]{TomanovWeiss} asked whether
it is possible for there to be a divergent $\T_{\real}$-orbit when $\dim
\T > \Qrank \G$. B.~Weiss \cite[Conj.~4.11A]{Weiss-DivTraj} conjectured
that the answer is negative.

\begin{conj} \label{DivConj}
 Let
 \begin{itemize}
 \item $\G$ be a semisimple algebraic group that is defined
over~$\rational$,
 \item $\Gamma$ be a subgroup of $\G_{\real}$ that is commensurable with
$\G_\integer$,
 \item $T$ be a connected Lie subgroup of an $\real$-split
torus in~$\G_{\real}$,
 and
 \item $x_0 \in \G_{\real}$.
 \end{itemize}
 If the $T$-orbit of~$\Gamma x_0$ is divergent in $\Gamma
\backslash \G_{\real}$, then $\dim T \le \Qrank \G$.
 \end{conj}

The conjecture easily reduces to the case where $\G$ is connected and
$\rational$-simple. Furthermore, the desired conclusion is obvious if
$\Qrank \G = 0$ (because this implies that $\Gamma \backslash \G_{\real}$
is compact), and it is easy to prove if $\Qrank \G = 1$
\seeSect{Qrank1PfSect}. Our main result is that the conjecture is also
true in the first interesting case:

\begin{thm} \label{TDivRank2}
 Suppose $\G$, $\Gamma$, $T$, and~$x_0$ are as specified in
Conj.~\ref{DivConj}, and assume $\Qrank \G \le 2$.  If the $T$-orbit
of~$\Gamma x_0$ is divergent in $\Gamma \backslash \G_{\real}$, then
$\dim T \le \Qrank \G$.
 \end{thm}

For higher $\Qrank$s,  we prove only the upper bound $\dim T < 2 (\Qrank
\G)$ \see{HighQrankThm}. A result of G.~Tomanov and B.~Weiss
\cite[Thm.~1.4]{TomanovWeiss} asserts that if $\Qrank \G < \Rrank \G$,
then $\dim T < \Rrank \G$. After seeing a preliminary version of our
work, B.~Weiss \cite{Weiss-Qrank} has recently proved the conjecture in
all cases.

\subsection*{Geometric reformulation}
 We remark that, by using the well-known fact that flats in a symmetric
space of noncompact type are orbits of $\real$-split tori in its isometry
group \cite[Prop.~6.1, p.~209]{Helgason}, the conjecture and our theorem
can also be stated in the following geometric terms.

Suppose $\cover{X}$ is a symmetric space, with no Euclidean (local)
factors. Recall that a \emph{flat} in~$\cover{X}$ is a connected, totally
geodesic, flat submanifold of~$\cover{X}$.  Up to isometry, $\cover{X} =
G/K$, where $K$ is a compact subgroup of a connected, semisimple Lie
group $G$ with finite center. Then $\Rrank G$ has the following
geometric interpretation:

\begin{fact} \label{Rrank-geometric}
 $\Rrank G$ is the largest natural number~$r$, such that $\cover{X}$
contains a topologically closed, simply connected, $r$-dimensional flat.
 \end{fact}

Now let $X = \Gamma \backslash \cover{X}$ be a locally symmetric space
modeled on~$X$, and assume that $X$ has finite volume. Then $\Qrank
\Gamma$ is a certain algebraically defined invariant of~$\Gamma$.
It can be characterized by
the following geometric property:

\begin{prop} \label{Qrank=BddDist}
 $\Qrank \Gamma$ is the smallest natural number~$r$, for which there
exists collection of finitely many $r$-dimensional flats in~$X$, such
that all of $X$ is within a bounded distance of the union of these flats. 
 \end{prop}

 It is clear from this that the $\Qrank$ does not change if $X$~is
replaced by a finite cover, and that it satisfies $\Qrank \Gamma \le
\Rrank G$. Furthermore, the algebraic definition easily implies that if
$\Qrank \Gamma = r$, then some finite cover of~$X$ contains a
topologically closed, simply connected flat of dimension~$r$. If
Conj.~\ref{DivConj} is true, then there are no such flats of larger
dimension. In other words, $\Qrank$ should have the following geometric
interpretation, analogous to \pref{Rrank-geometric}:

\begin{conj} \label{Qrank-geometric}
 $\Qrank \Gamma$ is the largest natural number~$r$, such that some finite
cover of~$X$ contains a topologically closed, simply connected,
$r$-dimensional flat.
 \end{conj}

More precisely, Conj~\ref{DivConj} is equivalent to the assertion that
$\Qrank \Gamma$ is the largest natural number~$r$, such that $\cover{X}$
contains a topologically closed, simply connected, $r$-dimensional
flat~$F$, for which the composition $F \hookrightarrow \cover{X} \to X$
is a proper map.

\begin{ack}
 The authors would like to thank Kevin Whyte for helpful discussions
related to Prop.~\ref{ComponentsProp}. D.~W.~M.\ was partially supported
by a grant from the National Science Foundation (DMS--0100438).
 \end{ack}

\section{Example: a proof for $\Qrank 1$} \label{Qrank1PfSect}

To illustrate the ideas in our proof of Thm.~\ref{TDivRank2}, we sketch a
simple proof that applies when $\Qrank \G = 1$. (A similar proof appears
in \cite[Prop.~4.12]{Weiss-DivTraj}.)

\begin{proof}
 Suppose $\G$, $\Gamma$, $T$, and~$x_0$ are as specified in
Conj.~\ref{DivConj}. For convenience, let $\pi \colon \G_\real \to \Gamma
\backslash \G_\real$ be the natural covering map. Assume that $\Qrank \G
= 1$, that $\dim T = 2$, and that the $T$-orbit of~$\pi(x_0)$ is
divergent in $\Gamma \backslash G$. This will lead to a contradiction.

Let $E_1 = \Gamma \backslash \G_{\real}$. Because $\Qrank \G = 1$,
reduction theory (the theory of Siegel sets) implies that there exist
 \begin{itemize}
 \item a compact subset~$E_0$ of $\Gamma \backslash \G_{\real}$, 
 and
 \item a $\rational$-representation $\rho \colon \G \to {\bf GL}_n$ (for
some~$n$),
 \end{itemize}
 such that, for each connected component~$\mathcal{E}$ of $G_\real
\smallsetminus \pi^{-1}(E_0)$, there is a nonzero vector $v \in
\rational^n$, such that
 \begin{equation} \label{Qrank1-to0Eqn}
 \mbox{if $\Gamma g_n \to \infty$ in $\Gamma \backslash \G_\real$, and
$\{g_n\} \subset \mathcal{E}$, then $\rho(g_n) v \to 0$.}
 \end{equation}
 (In geometric terms, this is the fact that, because $E_1 \smallsetminus
E_0$ consists of disjoint ``cusps," $\G_{\real} \smallsetminus
\pi^{-1}(E_0)$ consists of disjoint ``horoballs.")

Given $\epsilon > 0$, let $T_R$ be a large circle ($1$-sphere) in~$T$,
centered at the identity element. Because the $T$-orbit of~$\pi(x_0)$ is
divergent, we may assume $\pi(x_0 T_R)$ is disjoint from~$E_0$. Then,
because $T_R \approx S^1$ is connected, the set $x_0 T_R$ must be
contained in a single component of $G_\real \smallsetminus
\pi^{-1}(E_0)$. Thus, there is a vector $v \in \rational^n$, such that
$\|\rho(t) v\| < \epsilon \|v\|$ for all $t \in T_R$.

Fix some $t \in T_R$. Then $t^{-1}$ also belongs to~$T_R$, so $\|\rho(t)
v\|$ and $\|\rho(t^{-1}) v\|$ are both much smaller than~$\|v\|$. This is
impossible \see{AntipodContract}.
 \end{proof}

The above proof does not apply directly when $\Qrank \G = 2$, because, in
this case, there are arbitrarily large compact subsets~$C$ of~$\Gamma
\backslash \G_{\real}$, such that $\G_\real \smallsetminus \pi^{-1}(C)$ is
connected. Instead of only $E_0$ and~$E_1$, we consider a more refined
stratification $E_0 \subset E_1 \subset E_2$ of $\Gamma \backslash G$.
(It is provided by the structure of Siegel sets in $\Qrank$ two. The set
$E_0$ is compact, and, for $i \ge 1$, each component~$\mathcal{E}$ of
$\pi^{-1}(E_i \smallsetminus E_{i-1})$ has a corresponding
representation~$\rho$ and vector~$v$, such that \pref{Qrank1-to0Eqn}
holds. Thus, it suffices to find a component of either $\pi^{-1}(E_1
\smallsetminus E_0)$ or $\pi^{-1}(E_2 \smallsetminus E_1)$ that contains
two antipodal points of~$T_R$. Actually, we replace $E_1$ with a slightly
larger set that is open, so that we may apply the following property
of~$S^2$:

\begin{prop}[{see \ref{OpenSetsProp}}] \label{ComponentsProp}
 Suppose $n \ge 2$, and that $\{V_1,V_2\}$ is an open cover of the
$n$-sphere~$S^n$ that consists of only $2$~sets. Then there is a
connected component~$C$ of some~$V_i$, such that $C$ contains two
antipodal points of~$S^n$.
 \end{prop}

\begin{rem}
 In \S\ref{MainPfSect}, we do not use the notation $E_0 \subset E_1
\subset E_2$. The role of~$E_0$ is played by $\pi(Q \Siegel_\Delta^+)$, 
the role of an open set containing~$E_1$ is played by $\pi(Q
\Siegel_\alpha \cup Q \Siegel_\beta)$, and the role of $E_2
\smallsetminus E_1$ is played by $\pi(Q \Siegel_*)$.
 \end{rem}

\section{Preliminaries}

The classical Borsuk-Ulam Theorem implies that if $f \colon S^n \to
\real^k$ is a continuous map, and $n \ge k$, then there exist two
antipodal points $x$ and~$y$ of~$S^n$, such that $f(x) = f(y)$. We use
this to prove the following stronger version of
Prop.~\ref{ComponentsProp}:

 \begin{prop} \label{OpenSetsProp}
 Suppose $\opencover$ is an open cover of~$S^n$, with $n \ge 2$, such that
no point of~$S^n$ is contained in more than~$2$ of the sets
in~$\opencover$. Then some $V \in \opencover$ contains two antipodal
points of~$S^n$.
 \end{prop}

\begin{proof}
 Because $S^n$ is compact, we may assume $\opencover$ is finite.
 Let $\{\phi_V\}_{V \in \opencover}$ be a partition of unity subordinate
to~$\opencover$. This naturally defines a continuous function $\Phi$
from~$S^n$ to the simplex
 $$ \Delta_{\opencover} = \bigset{ (x_V)_{V \in \opencover} }{ \sum_{V
\in \opencover} x_V = 1 }
 \subset [0,1]^{\opencover} .$$
 Namely, 
 $\Phi(x) = \bigl( \phi_V(x) \bigr)_{V \in \opencover}$. 
 Our hypothesis on~$\opencover$ implies that no more than $2$~components
of $\Phi(x)$ are nonzero, so the image of~$\Phi$ is contained in the
$1$-skeleton $\Delta_{\opencover}^{(1)}$ of~$\Delta_{\opencover}$.
Because $S^n$ is simply connected, $\Phi$ lifts to a map from~$S^n$ to
the universal cover $\widetilde{\Delta_{\opencover}^{(1)}}$ of
$\Delta_{\opencover}^{(1)}$. The universal cover is a tree, which can be
embedded in~$\real^2$, so the Borsuk-Ulam Theorem implies that there
exist two antipodal points $x$ and~$y$ of~$S^n$, such that
$\widetilde\Phi(x) = \widetilde\Phi(y)$. Thus, there exists $V \in
\opencover$, such that $\phi_V(x) = \phi_V(y) \neq 0$. So $x,y \in V$.
 \end{proof}

For completeness, we also provide a proof of the following simple
observation.

\begin{lem} \label{AntipodContract}
  Let $T$ be any abelian group of diagonalizable $n \times n$ real
matrices.  There is a constant $\epsilon > 0$, such that
	if  
 \begin{itemize}
 \item $v$  is any vector in  $\real^n$,
	and 
 \item $t$  is any element of~$T$,
 \end{itemize}
	then either  $\|t v\| \ge \epsilon \|v\|$  or  $\|t^{-1} v\|  \ge
\epsilon \|v\|$.
 \end{lem}

\begin{proof}
  The elements of~$T$ can be simultaneously diagonalized.
  Thus, after a change of basis (which affects norms by only a bounded
factor), we may assume that each standard basis vector~$e_i$ is an
eigenvector for every element of~$T$.

	Let $\epsilon = 1/n$, write $v = (v_1,\ldots,v_n)$, and let $t_i$  be
the eigenvalue of~$t$  corresponding to the eigenvector~$e_i$.	Some
component~$v_j$ of~$v$ must be at least~$\|v\|/n$ in absolute value. We
may  assume $|t_j| \ge 1$, by replacing $t$  with~$t^{-1}$  if
necessary.  Then
	$$ \|t v\| = \|(t_1 v_1,\ldots,t_n v_n)\| \ge |t_j v_j| \ge 
1 \cdot \frac{\|v\|}{n} = \epsilon \|v\| ,$$
  as desired. 
 \end{proof}

\section{Properties of Siegel sets}

We present some basic results from reduction theory that follow easily
from the fundamental work of A.~Borel and Harish-Chandra
\cite{BorelHarishChandra} (see also \cite[\S13--\S15]{Borel-IntroArith}).
Most of what we need is essentially contained in
\cite[\S2]{Leuzinger-exhaust}, but we are working in~$G$, rather than in
$\cover{X} = G/K$. We begin by setting up the standard notation.

\begin{notation}[{cf.\ \cite[\S1]{Leuzinger-exhaust}}]
 Let
 \begin{itemize}
 \item $\G$ be a connected, almost simple $\rational$-group, with $\Qrank
\G = 2$,
 \item $G$ be the identity component of~$\G_{\real}$,
 \item $\Gamma$ be a finite-index subgroup of $G_\integer \cap G$,
 \item $P$ be a minimal parabolic $\rational$-subgroup of~$G$,
 \item $\A$ be a maximal $\rational$-split torus of~$\G$,
 \item  $A$ be the identity component of~$\A_{\real}$,
 and
 \item $K$ be a maximal compact subgroup of~$G$.
 \end{itemize}
 We may assume $A \subset P$. Then we have a Langlands decomposition $P =
UMA$, where $U$ is unipotent and $M$ is reductive. We remark that $U$
and~$A$ are connected, but $M$ is not connected (because $P$ is not
connected).
 \end{notation}

\begin{notation}[{cf.\ \cite[\S1]{Leuzinger-exhaust}}]
 The choice of~$P$ determines an ordering of the $\rational$-roots
of~$\G$. Because $\Qrank G = 2$, there are precisely two simple
$\rational$-roots $\alpha$ and~$\beta$ (so the base $\Delta$ is
$\{\alpha,\beta\}$). Then $\alpha$ and~$\beta$ are homomorphisms from $A$
to~$\real^+$. 

Any element~$g$ of~$G$ can be written in the form $g = pak$, with $p \in
UM$, $a \in A$, and $k \in K$. The element~$a$ is uniquely determined
by~$g$, so we may use this decomposition to extend $\alpha$ and~$\beta$
to continuous functions  $\tilde\alpha$ and $\tilde\beta$ defined on all
of~$G$:
 \begin{align*}
 \tilde\alpha(g) &= \alpha(a) \text{ \ if $g \in UMaK$ and $a \in A$}, \\
 \tilde\beta(g) &= \beta(a) \text{ \ if $g \in UMaK$ and $a \in A$}.
 \end{align*}
 \end{notation}

\begin{notation}[{cf.\ \cite[\S2]{Leuzinger-exhaust}}] \ 
 \begin{itemize}

 \item Fix a subset $Q$ of $\G_{\rational} \cap G$, such that
 $$ \text{$Q$ is a set of representatives of $\Gamma \backslash
(\G_{\rational} \cap G) / (\Pg_{\rational} \cap P)$.} $$
 Note that $Q$ is finite.

 \item For $\tau > 0$, let 
 $ A_\tau = \{\, a \in A \mid \text{$\alpha(a) > \tau$ and $\beta(a) >
\tau$} \,\}$.

 \item For $\tau > 0$ and a precompact, open subset $\omega$ of $UM$, let
$\Siegel_{\tau,\omega} = \omega A_\tau K$. This is a
\emph{Siegel set} in~$G$.

 \item We fix $\tau > 0$ and a precompact, open subset $\omega$ of $UM$,
such that, letting $\Siegel = \Siegel_{\tau,\omega}$, we have
 $$ \text{$Q \Siegel$ is a fundamental set for $\Gamma$ in~$G$} .$$
 That is, 
 \begin{itemize}
 \item $\Gamma Q \Siegel = G$,
 and 
 \item $\{\, \gamma \in \Gamma \mid \gamma Q \Siegel \cap pQ \Siegel \neq
\emptyset \,\}$ is finite, for all $p \in G_{\rational} \cap G$.
 \end{itemize}

 \item Let $\mathcal{D} = \bigset
 {p^{-1} \gamma q}
 {\begin{matrix}
 p,q \in Q, \ \gamma \in \Gamma,
 \\
 \text{$p \Siegel \cap \gamma q \Siegel$ is \emph{not} precompact}
 \end{matrix}
 } \subset \G_{\rational} \cap G$.
 Note that $\mathcal{D}$ is finite.

\item Fix $r > 0$, such that, for $q \in \mathcal{D}$, we have
 \begin{itemize}
 \item if $\tilde\alpha$ is bounded on $\Siegel \cap q
\Siegel$, then $\tilde\alpha(\Siegel \cap q
\Siegel) < r$,
 and
 \item if $\tilde\beta$ is bounded on $\Siegel \cap q
\Siegel$, then $\tilde\beta(\Siegel \cap q
\Siegel) < r$.
 \end{itemize}

\item Fix any $r^* > r$.

 \item Define 
 \begin{itemize}
 \item $\Siegel_* = \{\, x \in \Siegel \mid \text{$\tilde\alpha(x) > r$
and $\tilde\beta(x) > r$} \,\}$,
 \item $\Siegel_\alpha = \{\, x \in \Siegel \mid \tilde\alpha(x) < r^*
\,\}$,
 \item $\Siegel_\beta = \{\, x \in \Siegel \mid \tilde\beta(x) < r^*
\,\}$,
 \and
 \item $\Siegel_\Delta = \Siegel_\alpha \cap \Siegel_\beta$.
 \end{itemize}
 Note that $\{\Siegel_*,\Siegel_\alpha,\Siegel_\beta\}$ is an open cover
of~$\Siegel$ (whereas \cite[p.~398]{Leuzinger-exhaust} defines
$\{\Siegel_*,\Siegel_\alpha,\Siegel_\beta\}$ to be a partition
of~$\Siegel$, so not all sets are open).
We have
 $$ G = \Gamma Q \Siegel_* \cup \Gamma Q \Siegel_\alpha \cup \Gamma
Q \Siegel_\beta .$$

\item For $p,q \in Q$, let 
 \begin{itemize}
 \item $ \mathcal{D}^{p,q}_0 = \{\, \gamma \in \Gamma \mid
\mbox{$p \Siegel \cap \gamma q \Siegel$ is precompact and nonempty} \,\} $,
 \item $\mathcal{D}^{p,q}_\alpha = \{\, \gamma \in \Gamma \mid
\mbox{$p \Siegel_\alpha \cap \gamma q \Siegel_\alpha$ is precompact and
nonempty} \,\} $,
 \item $\mathcal{D}^{p,q}_\beta = \{\, \gamma \in \Gamma \mid
\mbox{$p \Siegel_\beta \cap \gamma q \Siegel_\beta$ is precompact and
nonempty} \,\} $,
 \item $\mathcal{D}^{p,q}_{\alpha,\beta} = \{\, \gamma \in \Gamma \mid
\mbox{$p \Siegel_\alpha \cap \gamma q \Siegel_\beta$ is precompact and
nonempty} \,\} $,
 \end{itemize} 
 and, using an overline to denote the closure of a set,
 \begin{align*}
 \Siegel_\Delta^+ =
 \bigcup_{\gamma \in \mathcal{D}^{p,q}_0}
 & (\closure{p \Siegel \cap \gamma
q \Siegel})
 \cup \bigcup_{\gamma \in \mathcal{D}^{p,q}_\alpha} (\closure{p \Siegel_\alpha \cap
\gamma q \Siegel_\alpha})
 \\&
 \cup \bigcup_{\gamma \in \mathcal{D}^{p,q}_\beta}
(\closure{p \Siegel_\beta \cap \gamma q \Siegel_\beta})
 \cup \bigcup_{\gamma \in \mathcal{D}^{p,q}_{\alpha,\beta}}
(\closure{p \Siegel_\alpha \cap \gamma q \Siegel_\beta}) .
 \end{align*}
 Note that
$\mathcal{D}^{p,q}_0$, $\mathcal{D}^{p,q}_\alpha$, $\mathcal{D}^{p,q}_\beta$,
and~$\mathcal{D}^{p,q}_0$ are finite (because $Q\Siegel$ is a
fundamental set), so $\Siegel_\Delta^+$ is compact. And $\Gamma
\Siegel_\Delta^+$ is closed.
 \item For $\Theta \subset \Delta$, we use $P_\Theta$ to denote the
corresponding standard parabolic $\rational$-subgroup of~$G$
corresponding to~$\Theta$. In particular, $P_{\emptyset} = P$ and
$P_\Delta = G$. There is a corresponding Langlands decomposition
$P_\Theta = U_\Theta M_\Theta A_\Theta$.
 \end{itemize}
 \end{notation}

\begin{lem} \label{gSaMeetsSb}
 For all $\gamma \in \Gamma$ and $p,q \in Q$, we have:
 \begin{enumerate}
 \item \label{gSaMeetsSb-cpct}
 $p\Siegel_\alpha \cap \gamma q \Siegel_\beta$ is precompact,
 and
 \item \label{gSaMeetsSb-inSDelta+}
 $p \Siegel_\alpha \cap \gamma q
\Siegel_\beta \subset \Siegel_\Delta^+$.
 \end{enumerate}
 \end{lem}

\begin{proof}
 It suffices to prove \pref{gSaMeetsSb-cpct}, for then
\pref{gSaMeetsSb-inSDelta+} is immediate from the definition
of~$\Siegel_\Delta^+$ (and~$\mathcal{D}^{p,q}_{\alpha,\beta}$).
 Thus, let us suppose that $p \Siegel_\alpha \cap \gamma q \Siegel_\beta$
is not precompact.
 This will lead to a contradiction.

Because $\tilde\alpha$ is bounded on~$\Siegel_\alpha$, but $\Siegel_\alpha
\cap p^{-1} \gamma q \Siegel_\beta$ is not precompact, we know that
$\tilde\beta$ is unbounded on $\Siegel_\alpha \cap
p^{-1} \gamma q \Siegel_\beta$ (and, hence, on
$\Siegel \cap p^{-1} \gamma q \Siegel$). Therefore, 
\cite[Prop.~2.3]{Leuzinger-exhaust} implies that
 $$ p^{-1} \gamma q \in P_\alpha.$$
 Similarly (replacing $\gamma$ with~$\gamma^{-1}$ and interchanging $p$
with~$q$ and $\alpha$ with~$\beta$), because $ \gamma^{-1} p
\Siegel_\alpha \cap q \Siegel_\beta = \gamma^{-1}(p \Siegel_\alpha \cap
\gamma q \Siegel_\beta)$ is not precompact, we see that
 $$ q^{-1} \gamma^{-1} p \in P_\beta .$$
 Noting that $q^{-1} \gamma^{-1} p = (p^{-1} \gamma q)^{-1}$, we conclude
that $p^{-1} \gamma q \in P_\alpha \cap P_\beta = P_\emptyset$, so
\cite[Lem.~2.4(i)]{Leuzinger-exhaust} tells us that $p = q$ and
$p^{-1} \gamma q \in UM$.
 Therefore 
 $$\tilde\alpha(\Siegel_\alpha \cap p^{-1} \gamma q \Siegel_\beta )
\subset \tilde\alpha(\Siegel_\alpha)$$
 and
 $$\tilde\beta(\Siegel_\alpha \cap p^{-1} \gamma q \Siegel_\beta) \subset
 \tilde\beta( p^{-1} \gamma q \Siegel_\beta)
 \subset \tilde\beta( UM \Siegel_\beta) 
 = \tilde\beta(\Siegel_\beta) $$
 are precompact. So $\Siegel_\alpha \cap p^{-1} \gamma q \Siegel_\beta$ is
precompact, which contradicts our assumption that $p \Siegel_\alpha \cap
\gamma q \Siegel_\beta$ is not precompact.
 \end{proof}

\begin{lem} \label{gS*meetsS*}
 If $\gamma \in \Gamma$ and $p,q \in Q$, such that $p \Siegel_* \cap
\gamma q \Siegel_* \not \subset \Siegel_\Delta^+$, then $p = q$ and
$p^{-1} \gamma q \in (UM)_\rational$.
 \end{lem}

\begin{proof}
 It suffices to show that both $\tilde\alpha$ and~$\tilde\beta$ are
unbounded on $\Siegel \cap p^{-1} \gamma q \Siegel$, for then the desired
conclusion is obtained from \cite[Prop.~2.3 and
Lem.~2.4(i)]{Leuzinger-exhaust}. Thus, let us suppose (without loss of
generality) that 
 $$ \mbox{$\tilde\alpha$ is bounded on $\Siegel \cap p^{-1} \gamma q
\Siegel$.} $$
 This will lead to a contradiction.

\setcounter{case}{0}

\begin{case}
 Assume $\tilde\beta$ is also bounded on $\Siegel \cap p^{-1} \gamma q
\Siegel$.
 \end{case}
 Then $p \Siegel \cap \gamma q \Siegel = p(\Siegel \cap p^{-1} \gamma q
\Siegel)$ is precompact, so, by definition, $p \Siegel \cap \gamma q \Siegel
\subset \Siegel_\Delta^+$. Therefore
 $$ p \Siegel_* \cap \gamma q \Siegel_* \subset p \Siegel \cap \gamma
q \Siegel \subset \Siegel_\Delta^+ .$$
 This contradicts the hypothesis of the lemma.

\begin{case}
 Assume $\tilde\beta$ is not bounded on $\Siegel \cap p^{-1} \gamma q
\Siegel$.
 \end{case}
 From \cite[Lem.~2.5]{Leuzinger-exhaust}, we see that $p \Siegel \cap
\gamma q \Siegel \subset p \Siegel_\alpha$.  Therefore
 $$ p \Siegel_* \cap \gamma q \Siegel_* \subset p \Siegel_* \cap
p \Siegel_\alpha = \emptyset \subset \Siegel_\Delta^+ .$$
 This contradicts the hypothesis of the lemma.
 \end{proof}

\begin{cor} \label{gS*comps}
 If $x$ and~$y$ are two points in the same connected component of $\Gamma Q \Siegel_* \smallsetminus \Gamma
\Siegel_\Delta^+$, then there exist $\gamma_0,\gamma \in \Gamma$ and $q
\in Q$, such that
 $x \in \gamma_0 q \Siegel_*$, 
 $y \in \gamma_0 \gamma q \Siegel_*$,
 and $q^{-1} \gamma q \in (U M)_{\rational}$.
 \end{cor}

\begin{lem} \label{gSaMeetsSa} \ 
 \begin{enumerate}
 \item \label{gSaMeetsSa-meet}
 If $\gamma \in \Gamma$ and $p,q \in Q$, such that $p \Siegel_\alpha \cap
\gamma q \Siegel_\alpha \not \subset \Siegel_\Delta^+$, then $p^{-1}
\gamma q \in (P_\alpha)_\rational$.
 \item \label{gSaMeetsSa-inP}
 For each $p,q \in Q$, there exists $h_{p,q} \in
(P_{\alpha})_{\rational}$, such that $p^{-1} \Gamma q \cap
(P_{\alpha})_{\rational} \subset h_{p,q} (U_\alpha
M_\alpha)_\rational$. 
 \end{enumerate}
 \end{lem}

\begin{proof}
 \pref{gSaMeetsSa-meet} Because $p \Siegel_\alpha \cap \gamma q
\Siegel_\alpha \not \subset \Siegel_\Delta^+$, we know, from the
definition of~$\Siegel_\Delta^+$ (and~$\mathcal{D}^{p,q}_\alpha$) that $p
\Siegel_\alpha \cap \gamma q \Siegel_\alpha$ is not precompact. Since
$\tilde\alpha$ is bounded on~$\Siegel_\alpha$, we conclude that
$\tilde\beta$ is \textbf{not} bounded on $\Siegel_\alpha \cap p^{-1}
\gamma q \Siegel_\alpha$ (and, hence, on $\Siegel \cap p^{-1} \gamma q
\Siegel$). Then \cite[Prop.~2.3]{Leuzinger-exhaust} asserts that $p^{-1}
\gamma q \in (P_\Theta)_\rational$, for $\Theta = \{\alpha\}$
or~$\emptyset$. Because $P_\emptyset \subset P_\alpha$, we conclude that
$p^{-1} \gamma q \in (P_\alpha)_\rational$.

 \pref{gSaMeetsSa-inP} From \cite[Lem.~2.4(ii)]{Leuzinger-exhaust}, we
see that the coset $(p^{-1} \gamma q) (U_\alpha M_\alpha)_\rational$ does
not depend on the choice of~$\gamma$, if we require $\gamma$ to be an
element of~$\Gamma$, such that $p^{-1} \gamma q \in
(P_{\alpha})_{\rational}$.
 \end{proof}

\begin{cor} \label{gSalphacomps}
 If $x$ and~$y$ are two points in the same connected component of $\Gamma
Q \Siegel_\alpha \smallsetminus \Gamma \Siegel_\Delta^+$, then there exist
$\gamma_0,\gamma \in \Gamma$ and $p,q \in Q$, such that
 $x \in \gamma_0 p \Siegel_\alpha$, 
 $y \in \gamma_0 \gamma q \Siegel_\alpha$,
 and $p^{-1} \gamma q \in h_{p,q}(U_\alpha M_\alpha)_{\rational}$.
 \end{cor}

\section{Proof of the Main Theorem} \label{MainPfSect}

Let $G$, $\Gamma$, $T$ and~$x_0$ be as described in the
hypotheses of Thm.~\ref{TDivRank2}, and assume $\dim T \ge 3$. (This will
lead to a contradiction.) Let $T_R$ be a large sphere (centered at the
origin) in~$T$. Because $\Siegel_\Delta^+$ is compact and the $T$-orbit
of~$\Gamma x_0$ is divergent in $\Gamma \backslash G$, we may assume that 
 \begin{equation}
 (x_0 T_R) \cap (\Gamma \Siegel_\Delta^+) = \emptyset .
 \end{equation}
 Let 
 $$ V_1 = \{\, t \in T_R \mid x_0 t \in \Gamma Q \Siegel_* \,\} $$
 and
 $$ V_2 = \{\, t \in T_R \mid x_0 t \in \Gamma Q \Siegel_\alpha \cup
\Gamma Q \Siegel_\beta \,\} . $$
 From Prop.~\ref{ComponentsProp}, we know there exists $t \in T_R$, and
a connected component~$C$ of either $V_1$ or $V_2$, such
that $t$ and~$t^{-1}$ both belong to~$C$.

\setcounter{case}{0}

\begin{case} \label{TDivRank2Case-notE}
 Assume $C$ is a component of $V_1$.
 \end{case}
 From Cor.~\ref{gS*comps}, we see that there exist $\gamma_0,\gamma \in
\Gamma$ and $q \in Q$, such that
 $x_0 t \in \gamma_0 q \Siegel_*$, 
 $x_0 t^{-1} \in \gamma_0 \gamma q \Siegel_*$,
 and $q^{-1} \gamma q \in (U M)_{\rational}$.
 Because $\Gamma x_0 t$ and $\Gamma x_0 t^{-1}$ are near infinity in
$\Gamma \backslash G$, we must have 
 \begin{enumerate}
 \item \label{mainpf-alpha(t)} either $\tilde\alpha(q^{-1} \gamma_0^{-1}
x_0 t) \gg 1$ or $\tilde\beta(q^{-1} \gamma_0^{-1} x_0 t) \gg 1$, and
 \item \label{mainpf-alpha(tmin)}
 either $\tilde\alpha(q^{-1} \gamma^{-1} \gamma_0^{-1} x_0 t^{-1}) \gg 1$
or $\tilde\beta(q^{-1} \gamma^{-1} \gamma_0^{-1} x_0 t^{-1}) \gg 1$.
 \end{enumerate}
 Since $q^{-1} \gamma q \in (U M)_{\rational}$ is sent to
the identity element by both~$\tilde\alpha$ and~$\tilde\beta$, we have
 \begin{enumerate}
 \item[(\ref{mainpf-alpha(tmin)}$'$)] \label{mainpf-alpha(tmin)'}
 either $\tilde\alpha(q^{-1} \gamma_0^{-1}
x_0 t^{-1}) \gg 1$ or $\tilde\beta(q^{-1} \gamma_0^{-1} x_0 t^{-1}) \gg
1$.
 \end{enumerate}

 Let  
 \begin{itemize}
 \item $V = \bigwedge^d \Lie G$, where $d = \dim U$,
 \item $\rho \colon G \to \GL(V)$ be the natural adjoint representation
of~$G$ on~$V$,
 \item $v_{\Lie U}$ be a nonzero element of~$V_\integer$ in the
one-dimensional subspace $\bigwedge^d \Lie U$,
 and
 \item $v_{\Lie U}' = \rho(x_0^{-1} \gamma_0 q) v_{\Lie U}$.
 \end{itemize}
 It is important to note that $\|v_{\Lie U}'\|$ is bounded away from $0$,
independent of the choice of~$q$ and~$\gamma_0$. (There are only
finitely many choices of~$q$, so $q$~is not really an issue. The key
point is that $\rho(q) v_{\Lie U}$ is a $\rational$-element of~$V$, so its
$\G_{\integer}$-orbit is bounded away from~$0$.)

 On the other hand, for any $g \in P_{\emptyset}$, we have $\rho(g^{-1})
v_{\Lie U} = \tilde\alpha(g)^{-\ell_1} \tilde\beta(g)^{-\ell_2} v_{\Lie
U}$, for some positive integers $\ell_1$ and~$\ell_2$ (because the sum of
the positive $\rational$-roots of~$\G$ is $\ell_1 \alpha + \ell_2 \beta$).
Therefore, from \pref{mainpf-alpha(t)} and (\ref{mainpf-alpha(tmin)}$'$),
we see that
 $$ \rho(t^{-1}) v_{\Lie U}' = \rho \bigl( (q^{-1} \gamma_0^{-1} x_0 t)^{-1}
\bigr) v_{\Lie U} \approx 0 $$
 and
 $$ \rho(t) v_{\Lie U}' = \rho \bigl( (q^{-1} \gamma_0^{-1} x_0 t^{-1})^{-1}
\bigr) v_{\Lie U} \approx 0 .$$
 This contradicts Lem.~\ref{AntipodContract}.

\begin{case}
 Assume $C$ is a component of~$V_2$.
 \end{case}
 From Lem.~\fullref{gSaMeetsSb}{inSDelta+}, we see that $x_0 C$ is
contained in either $\Gamma Q \Siegel_\alpha$ or  $\Gamma Q
\Siegel_\beta$. Assume, without loss of generality, that $x_0 C \subset
\Gamma Q \Siegel_\alpha$.
  From Cor.~\ref{gSalphacomps}, we see that there
exist $\gamma_0,\gamma \in \Gamma$ and $p,q \in Q$, such that
 $$\mbox{
 $x_0 t \in \gamma_0 p \Siegel_\alpha$, 
 $x_0 t^{-1} \in \gamma_0 \gamma q \Siegel_\alpha$,
 and $p^{-1} \gamma q \in h_{p,q} (U_\alpha M_\alpha)_{\rational}$.
 }$$
 Let $\Lie U_\alpha$ be the Lie algebra of~$U_\alpha$, and let
$\rho_\alpha \colon G \to \GL(V_\alpha)$ be the natural adjoint
representation of~$G$ on~$V_\alpha = \bigwedge^{d_\alpha} \Lie G$, where
$d_\alpha = \dim \Lie U_\alpha$.

We can obtain a contradiction by arguing as in
Case~\ref{TDivRank2Case-notE}, with the representation~$\rho_\alpha$ in
the place of~$\rho$. To see this, note that:
 \begin{itemize}

\item For $a \in \ker \alpha$, we have $\rho_\alpha(a^{-1}) v_{\Lie
U_\alpha} = \beta(a)^{-\ell} v_{\Lie U_\alpha}$, for some positive
integer~$\ell$. Since $\rho_\alpha(UM) \subset \rho_\alpha(U_\alpha
M_\alpha)$ fixes~$v_{\Lie U_\alpha}$, and $\rho_\alpha(K)$ is compact,
this implies that
 $$ \mbox{$\| \rho_\alpha(g^{-1}) v_{\Lie U_\alpha}\| \approx
\tilde\beta(g)^{-\ell} \|v_{\Lie U_\alpha}\|$
 \qquad for $g \in \Siegel_\alpha$.} $$

 \item Because $\Gamma x_0 t$ and $\Gamma x_0 t^{-1}$ are near infinity
in $\Gamma \backslash G$, and $\tilde\alpha$ is bounded on
$\Siegel_\alpha$, we must have 
 $$ \mbox{\pref{mainpf-alpha(t)}
 $\tilde\beta(p^{-1} \gamma_0^{-1} x_0 t) \gg 1$
 \qquad and \qquad
 \pref{mainpf-alpha(tmin)}
 $\tilde\beta(q^{-1} \gamma^{-1} \gamma_0^{-1} x_0 t^{-1}) \gg 1$.}$$
 Therefore, letting $v_{\Lie U_\alpha}' = \rho_\alpha(x_0^{-1} \gamma_0 p)
v_{\Lie U_\alpha}$, we have
 $$\text{(\ref{mainpf-alpha(t)}*) }
 \rho_\alpha(t^{-1}) v_{\Lie U_\alpha}' \approx 0. $$
 Because $h_{p,q} \in P_\alpha$ normalizes~$U_\alpha$, we have
$\rho_\alpha (h_{p,q} v_{\Lie U_\alpha}) = c_{p,q} v_{\Lie U_\alpha}$,
for some scalar~$c_{p,q}$. 
 Therefore, since $(p^{-1} \gamma q) h_{p,q}^{-1} \in (U_\alpha
M_\alpha)_{\rational}$ fixes $v_{\Lie U_\alpha}$, and $\{c_{p,q}\}$,
being finite, is bounded away from~$0$, we have
   $$
  \text{(\ref{mainpf-alpha(tmin)}*) }
  \rho_\alpha(t) v_{\Lie U_\alpha}'
    = \rho_\alpha\bigl( t x_0^{-1} \gamma_0 p) v_{\Lie U_\alpha}
  =  \rho_\alpha\bigl( t x_0^{-1} \gamma_0 p (p^{-1} \gamma
q) {h_{p,q}}^{-1} \bigr) v_{\Lie U_\alpha}
 =c_{p,q}^{-1} \, \rho_\alpha\bigl( t x_0^{-1}
\gamma_0 \gamma q) v_{\Lie U_\alpha}
 \approx 0 . 
 $$

 \end{itemize}

This completes the proof of Thm.~\ref{TDivRank2}.

\section{Results for higher $\Qrank$} \label{HighQrankSect}

The proof of Thm.~\ref{TDivRank2} generalizes to establish the following
result:

\begin{thm} \label{HighQrankThm}
 Suppose $\G$, $\Gamma$, $T$, and~$x_0$ are as specified in
Conj.~\ref{DivConj}, and assume $\Qrank \G \ge 1$. 
 If the $T$-orbit of~$\Gamma x_0$ is divergent in $\Gamma \backslash
\G_{\real}$, then $\dim T \le 2(\Qrank \G) - 1$.
 \end{thm}

\begin{proof}[Sketch of proof] 
 As in \cite[\S1 and \S2]{Leuzinger-exhaust}, let $\Delta$ be the set of
simple $\rational$-roots, construct a fundamental set $Q \Siegel$, define
the finite set~$\mathcal{D}$, and choose $r > 0$, such that, for $q \in
\mathcal{D}$ and $\alpha \in \Delta$, we have
 $$ \text{if $\tilde\alpha$ is bounded on $\Siegel \cap q \Siegel$, then
$\tilde\alpha(\Siegel \cap q \Siegel) < r$.} $$
 
Fix any $r^* > r$.
 For each subset $\Theta$ of~$\Delta$, let 
 $$ \text{$\Siegel_\Theta =  \{\, x \in \Siegel \mid \tilde\alpha(x) <
r^*, \ \forall \alpha \in \Theta \,\} $
 and
  $ \Siegel_\Theta^- =  \{\, x \in \Siegel \mid \tilde\alpha(x) \le r, \
\forall \alpha \in \Theta \,\}$,} $$
 and choose $h_{p,q}^\Theta$ such that $p^{-1} \Gamma q \cap
(P_\Theta)_\rational \subset h_{p,q}^\Theta (U_\Theta
M_\Theta)_\rational$ for $p,q \in Q$.
  Set $d = \Qrank \G$, and, for $i = 0,\ldots,d$, let
 $$ \text{$\displaystyle E_i = \bigcup_{\begin{matrix} \Theta \subset
\Delta \\ \#\Theta = i \end{matrix}} \Siegel_\Theta$
 \qquad and \qquad
 $\displaystyle E_i^- = \bigcup_{\begin{matrix} \Theta \subset \Delta \\
\#\Theta = i \end{matrix}} \Siegel_\Theta^-$.} $$
 Then $\{\, E_1, E_2 \smallsetminus E_1^-, \cdots, E_d \smallsetminus
E_{d-1}^-,\}$ is an open cover of $\Gamma \backslash G$.

 For $p,q \in Q$ and $\Theta_1,\Theta_2 \subset \Delta$, let
 $$ \mathcal{D}^{p,q}_{\Theta_1,\Theta_2} = \{\, \gamma \in \Gamma \mid
\mbox{$p \Siegel_{\Theta_1} \cap \gamma q \Siegel_{\Theta_2}$ is
precompact and nonempty}\,\} .$$
 Define 
 $$ \Siegel_\Delta^+ = \bigcup_{\begin{matrix}
 p,q \in Q, \\
 \Theta_1,\Theta_2 \subset \Delta, \\
 \gamma \in \mathcal{D}^{p,q}_{\Theta_1,\Theta_2}
 \end{matrix}}
 (\closure{p \Siegel_{\Theta_1} \cap \gamma q \Siegel_{\Theta_2}}) .$$

Suppose $\dim T \ge 2d$. Then we may choose a large $(2d-1)$-sphere $T_R$
in~$T$. Prop.~\ref{MoreComponentsProp} below implies that there exists $t
\in T_R$ and a component~$C$ of some $E_i \smallsetminus E_{i-1}^-$, such
that $x_0 t$ and~$x_0 t^{-1}$ belong to~$C$. If $T_R$ is chosen large
enough that $x_0 T_R$ is disjoint from $\Gamma \Siegel_\Delta^+$, then
there exist $\Theta \subset \Delta$ (with $\#\Theta = i$),
$\gamma_0,\gamma \in \Gamma$, and $p,q \in Q$, such that
 $x_0 t \in \gamma_0 p \Siegel_\Theta$, 
 $x_0 t^{-1} \in \gamma_0 \gamma q \Siegel_\Theta$,
 and $p^{-1} \gamma q \in h_{p,q}^\Theta (U_\Theta M_\Theta)_{\rational}$.
We obtain a contradiction as in Case~\ref{TDivRank2Case-notE} of
\S\ref{MainPfSect}, using $\Lie U_\Theta$ in the place of~$\Lie U$.
 \end{proof}

The following result is obtained from the proof of
Prop.~\ref{OpenSetsProp}, by using the fact that any simplicial complex
of dimension $d-1$ can be embedded in $\real^{2d-1}$.

\begin{prop} \label{MoreComponentsProp}
 Suppose $n \ge 2d-1$, and that $\{V_1,V_2, \ldots,V_d\}$ is an open cover
of the $n$-sphere~$S^n$ that consists of only $d$~sets. Then there is a
connected component~$C$ of some~$V_i$, such that $C$ contains two
antipodal points of~$S^n$.
 \end{prop}

\begin{rem}
 For $k \ge 1$, it is known \cite{Scepin,IzydorekJaworowski} that there
exist a simplicial complex $\Sigma^k$ of dimension~$k$ and a continuous
map $f \colon S^{2k-1} \to \Sigma^k$, such that no two antipodal points
of~$S^{2k-1}$ map to the same point of~$\Sigma^k$. This implies that the
constant $2d-1$ in Prop.~\ref{MoreComponentsProp} cannot be improved to
$2d-3$.
 \end{rem}

\begin{rem}
 If $\Qrank G = 2$, then the conclusion of Thm.~\ref{TDivRank2} is
stronger than that of Thm.~\ref{HighQrankThm}. The improved bound in
\pref{TDivRank2} results from the fact that if $d = 2$, then the
universal cover of any $(d-1)$-dimensional simplicial complex embeds in
$\real^2 = \real^{2d-2}$. (See the proof of Prop.~\ref{OpenSetsProp}.)
When $d > 2$, there are examples of (simply connected)
$(d-1)$-dimensional simplicial complexes that embed only in
$\real^{2d-1}$, not~$\real^{2d-2}$.
 \end{rem}


\begin{thebibliography}{TW}

\bibitem[B]{Borel-IntroArith}
 A.~Borel:
 \emph{Introduction aux groupes arithm\'etiques},
 Actualit\'es Scientifiques et Industrielles, No. 1341,
 Hermann, Paris, 1969.


\bibitem[BH]{BorelHarishChandra}
 A.~Borel and Harish-Chandra:
 Arithmetic subgroups of algebraic groups.  
 \emph{Ann. of Math.} (2)  75 (1962) 485--535.

\bibitem[H]{Helgason}
 S.~Helgason:
 \emph{Differential Geometry and Symmetric Spaces},
 Academic Press, New York, 1962.

\bibitem[IJ]{IzydorekJaworowski}
 M.~Izydorek and J.~Jaworowski:
 Antipodal coincidence for maps of spheres into complexes,
 \emph{Proc. Amer. Math. Soc.} 123 (1995) 1947--1950.

\bibitem[L]{Leuzinger-exhaust}
 E.~Leuzinger:
 An exhaustion of locally symmetric spaces by compact submanifolds with
corners.
 \emph{Invent. Math.} 121 (1995) 389--410.

\bibitem[S]{Scepin}
 E.~V.~\v S\v cepin:
 On a problem of L.~A.~Tumarkin,
 \emph{Soviet Math. Dokl.} 15 (1974) 1024--1026.

\bibitem[TW]{TomanovWeiss}
 G.~Tomanov and B.~Weiss:
 Closed orbits for actions of maximal tori on homogeneous
spaces, \emph{Duke Math. J.} 119 (2003) 367--392.

\bibitem[W1]{Weiss-DivTraj}
 B.~Weiss:
 Divergent trajectories on noncompact parameter spaces,
 \emph{Geom. Funct. Anal.} 14  (2004) 94--149.

\bibitem[W2]{Weiss-Qrank}
 B.~Weiss:
 Divergent trajectories and $\Qrank$
 (preprint).


\end{thebibliography}
\end{document}